\documentstyle[12pt]{article}
\newcommand{\const}{\mathop{\rm const}\limits}

\newcommand{\diam}{\mathop{\rm diam}\limits}

\newcommand{\radii}{\mathop{\rm radii}\limits}

\newcommand{\Law}{\mathop{\rm Law}\limits}

\newcommand{\Var}{\mathop{\rm Var}\limits}

\begin{document}

\begin{center}

{\bf CENTRAL LIMIT THEOREM AND} \par

\vspace{4mm}

{\bf EXPONENTIAL TAIL ESTIMATIONS IN}\par

\vspace{4mm}

{\bf HYBRID  LEBESGUE-CONTINUOUS SPACES. }\par

\vspace{4mm}

 $ {\bf E.Ostrovsky^a, \ \ L.Sirota^b } $ \\

\vspace{4mm}

$ ^a $ Corresponding Author. Department of Mathematics and computer science, Bar-Ilan University, 84105, Ramat Gan, Israel.\\
\end{center}
E - mail: \ galo@list.ru \  eugostrovsky@list.ru\\
\begin{center}
$ ^b $  Department of Mathematics and computer science. Bar-Ilan University,
84105, Ramat Gan, Israel.\\

E - mail: \ sirota3@bezeqint.net\\

\vspace{3mm}
                    {\sc Abstract.}\\

 \end{center}

 \vspace{4mm}

 We  study the Central Limit Theorem (CLT) in the so-called  hybrid Lebesgue-continuous spaces
and tail behavior of normed sums of centered random independent variables (vectors) with values in these spaces. \par

  \vspace{4mm}

{\it Key words and phrases:} Central Limit Theorem (CLT), random vectors, variables (r.v.), processes and fields (r.f.);
Banach spaces, separability, law of distribution, Gaussian distribution, weak convergence of measures, Pisier's estimate,
hybrid   Lebesgue-continuous and continuous - Lebesque norms and spaces, norms, metric entropy, characteristical
functional, Rosenthal constants and inequalities, exponential  upper tail estimates, triangle (Minkowsky) inequality,
permutation inequalities, martingale, mixingale, stationarity.\\

\vspace{4mm}

{\it 2000 Mathematics Subject Classification. Primary 37B30, 33K55; Secondary 34A34,
65M20, 42B25.} \par

\vspace{4mm}

\section{Notations. Statement of problem.}

\vspace{3mm}

{\bf 1.}  Let  $  (B, ||\cdot||B )  $  be separable Banach space and  $ \{  \xi_j   \}, \xi = \xi_1, \ j = 1,2,\ldots  $
be a sequence of centered in the weak sense: $ {\bf E} (\xi_i,b) = 0 \ \forall b \in B^* $  of independent identical distributed
(i.; i.d.)  random variables (r.v.) (or equally random vectors, with at the same abbreviation r.v.)  defined on some non-trivial
probability  space $  (\Omega = \{\omega\}, F, {\bf P})   $ with values  in the space  $ B. $  Denote

 $$
 S(n) = n^{-1/2} \sum_{j=1}^n \xi_j, \ n = 1,2,\ldots. \eqno(1.1)
 $$

 If we suppose that the r.v. $ \xi  $ has a weak second moment:

 $$
 \forall b \in B^* \ \Rightarrow  (Rb,b) := {\bf E } (\xi,b)^2 < \infty,  \eqno(1.2)
 $$
then the characteristical functional (more exactly, the sequence of  characteristical functionals)

$$
\phi_{S(n)}(b) := {\bf E} e^{i \ (S(n), b) } \eqno(1.3)
$$
of $ S(n) $ converges as $ n \to \infty $ to the  characteristical functional of (weak, in general case)
Gaussian r.v. $ S = S(\infty)  $ with parameters $ (0,R):   $

$$
\lim_{n \to \infty} \phi_{S(n)}(b) = e^{ - 0,5 (Rb,b) }.
$$
  Symbolically:  $ S \sim N(0,R)  $ or $ \Law(S) = N(0,R).  $ The operator $ R = R_{S} $ is called
 the covariation operator, or variance of the r.v. $  S: $

$$
R = \Var(S);
$$
note that $  R = \Var(\xi). $\par

\vspace{4mm}

 We recall the classical definition of the CLT in the space $  B. $ \par

 \vspace{3mm}

 {\bf Definition 1.1.} {\it  We will say that the sequence $  \{ \xi_i  \} $  satisfies the CLT in the space $  B, $
 write:  $ \{\xi_j \} \in CLT = CLT(B) $ or   simple:  $ \xi \in CLT(B),  $  if the
 limit Gaussian r.v. $  S $  belongs to the  space $  B $ with probability one:  $  {\bf P} (S \in B) = 1 $ and the sequence
 of distributions $\Law(S(n))  $  converges weakly as $ n \to \infty $  to the distribution of the r.v. $ S = S(\infty):  $  }\par

$$
\lim_{n \to \infty} \Law(S(n)) = \Law(S). \eqno(1.4)
$$
 The  equality (1.4) imply that for any  continuous functional $ F: B \to R $

$$
\lim_{n \to \infty} {\bf  P} ( F(S(n)) < x  ) = {\bf P} (F(S) < x) \eqno(1.5)
$$
almost everywhere. \par

 In particular,

$$
\lim_{n \to \infty} {\bf  P} ( ||S(n)||B < x  ) = {\bf P} (||S||B < x), \ x > 0.
$$

 \vspace{3mm}

{\bf 2.} The problem of describing of necessary (sufficient) conditions  for the infinite - dimensional CLT in Banach space $  B $
has a long history; see, for instance, the  monographs \cite{Araujo1} - \cite{Ostrovsky1} and articles
\cite{Garling1} - \cite{Zinn1} ; see also reference therein.\par
 The applications of considered theorem in statistics and method Monte-Carlo  see, e.g. in
\cite{Frolov1} -  \cite{Ostrovsky403}.\par

\vspace{3mm}

{\bf 3.} The cornerstone of  this problem is to establish the {\it weak compactness } of the distributions generated
in the space $  B  $ by the sequence $  \{  S(n) \}: $

$$
\nu_n(D) = {\bf P} ( S(n) \in D),
$$
where $  D  $ is Borelian set in $  B; $ see \cite{Prokhorov1}; \cite{Billingsley1},  \cite{Billingsley2}. \par

\vspace{3mm}

{\bf 4.}  Let $ T = \{t \} $ be  precompact topological  metrizable space.  A concrete choice of the distance on the set $  T $
will be clarified below. \par

\vspace{3mm}

  Let $ (X = \{x\}, A,\mu) $ also be measurable spaces with sigma-finite {\it separable}
non - trivial measures $ \mu. $ The separability  denotes that  the metric space
$ A_k  $  relative the distance

$$
\rho(D_1, D_2) = \mu(D_1 \Delta D_2) = \mu(D_1 \setminus D_2) + \mu(D_2 \setminus D_1) \eqno(1.6)
$$
is separable.\par

 We denote as ordinary

 $$
 |f|_p = |f|L(p) = \left[ \int_X |f(x)|^p  \ \mu(dx) \right]^{1/p}, \ 1 \le p < \infty.
 $$
In particular, for the random variable $ \xi: \Omega \to R $

$$
|\xi|_p  = |\xi|L(p) = \left[ {\bf E} |\xi|^p   \right]^{1/p}.
$$

\vspace{5mm}

{\bf Definitions of hybrid Lebesgue - continuous spaces.} \par

\vspace{3mm}

 We will distinguish two types of such a spaces. \\

\vspace{3mm}

 {\bf Definition of continuous-Lebesgue (Lebesgue-Riesz) space}  $ CL(p) = C(T, L_p(X)).  $  \\

 \vspace{3mm}

We will say that the (measurable) function of two variables $ f = f(x,t), \ x \in X, \ t \in T $  belongs to the space
$ CL(p) = C(T, L_p(X)),$ if the map $ t \to f(\cdot,t), \ t \in T  $ is continuous in  the $ C(T) $ sense: \par

$$
\lim_{\epsilon \to 0+} \sup_{d(t,s) < \epsilon } \left[ \int_X |f(x,t) - f(x,s)|^p \ \mu(dx) \right]^{1/p} = 0.
$$

 The norm of the function $ f(\cdot,\cdot) $ in this space is defined as follows:

 $$
 ||f(\cdot, \cdot)|| C(T, L_p(X)) = ||f(\cdot, \cdot)|| CL(p) = \sup_{t \in T} |f(t, \cdot)|_p.
 $$

\vspace{3mm}

 {\bf Definition of  Lebesgue-continuous space }  $ L(p)C =  L_p(X)C.  $  \\

 \vspace{3mm}

We will say that the  bi-measurable function of two variables $ f = f(x,t), \ x \in X, \ t \in T $  belongs to the space
$ L(p)C =  L_p(X)C,$ if the following  norm is finite:

\vspace{3mm}

$$
||f(\cdot, \cdot)|| L_pC  = ||f(\cdot, \cdot)|| L_p(X)C = | \sup_{t \in T} |f(t,\cdot)|_p
$$
 and furthermore

$$
\lim_{\epsilon \to 0+} \left[ \int_X \sup_{d(t,s) < \epsilon } |f(x,t) - f(x,s)|^p \ \mu(dx) \right]^{1/p} = 0.
$$

 These spaces  are complete separable Banach function spaces. The detail investigation of these spaces see, e.g. in
a monograph \cite{Kufner1}, p. 113 - 119.\par
  They are used, for instance, in the theory of non-linear evolution Partial Differential  Equations, see
\cite{Fujita1},  \cite{Kato1},  \cite{Kato2}, \cite{Lions1}, \cite{Lions2}, \cite{Taylor1}. \par

\vspace{4mm}

{\bf 6. Our goal in this short article is  to obtain some sufficient conditions for tail estimation  for
normed sums of  random vectors and for the Central Limit theorem  in described below so-called Lebesgue-continuous
Banach spaces.}\\

\vspace{4mm}

{\bf 7.} These spaces are continuous extremal cases of the so-called mixed (anisotropic) Lebesgue - Riesz spaces.
 They was introduced and investigated in an article of Benedek A. and Panzone  R.  \cite{Benedek1};  see also
an article of  R.A. Adams  \cite{Adams1} and a classical monograph written by O.V.Besov, V.P.Il'in and S.M.Nikolskii
 \cite{Besov1}, chapters 1,2. \par

 In detail: let  $ p = (p_1, p_2, . . . , p_l) $ be $ l- $ dimensional vector such that $ 1 \le p_j < \infty.$ \par
Let also $ {X_k, A_k, \mu_k} $ be measurable spaces with sigma - finite separable measures $ \ \mu_k. $\\

 Recall that the anisotropic  (mixed) Lebesgue - Riesz space $ L_{ \vec{p}} $ consists on all the  total measurable
real valued function  $ f = f(x_1,x_2,\ldots, x_l) = f( \vec{x} ) $

$$
f:  \otimes_{k=1}^l X_k \to R
$$

with finite norm $ |f|_{ \vec{p} } \stackrel{def}{=} $

$$
\left( \int_{X_l} \mu_l(dx_l) \left( \int_{X_{l-1}} \mu_{l-1}(dx_{l-1}) \ldots \left( \int_{X_1}
 |f(\vec{x})|^{p_1} \mu(dx_1) \right)^{p_2/p_1 }  \ \right)^{p_3/p_2} \ldots   \right)^{1/p_l}.
$$

 Note that in general case $ |f|_{p_1,p_2} \ne |f|_{p_2,p_1}, $
but $ |f|_{p,p} = |f|_p. $ \par

 Observe also that if $ f(x_1, x_2) = g_1(x_1) \cdot g_2(x_2) $ (condition of factorization), then
$ |f|_{p_1,p_2} = |g_1|_{p_1} \cdot |g_2|_{p_2}, $ (formula of factorization). \par

 Note that under conditions separability of measures $ \{  \mu_k \} $   this spaces are also  separable and Banach spaces. \par

  These spaces arises in the Theory of Approximation, Functional Analysis, theory of Partial Differential Equations,
theory of Random Processes etc. \par

\vspace{3mm}

 Let for example $  l = 2; $ we agree to rewrite for clarity the expression for $ |f|_{p_1, p_2}  $ as follows:

$$
|f|_{p_1, p_2} := |  f|_{p_1, X_1; p_2, X_2}.
$$
 Analogously,

$$
|f|_{p_1, p_2,p_3} = |  f|_{p_1, X_1; p_2, X_2; p_3, X_3}.
$$

 Note that under imposed condition of continuity

$$
||f(\cdot, \cdot)|| C(T, L_p(X)) =  |f|_{p,X; T, \infty}
$$
and

$$
||f(\cdot, \cdot)|| L_p(X)C = |f|_{T, \infty; p,X}.
$$

  It is known \cite{Adams1}  see also \cite{Besov1}, chapter 1,  p. 24 - 26, ("permutation inequality"), that

$$
||f(\cdot, \cdot)|| C(T, L_p(X)) \le  |f|_{T, \infty; p,X}.
$$

 The CLT in mixed $ L_{\vec{p}}  $ spaces  is considered in \cite{Ostrovsky303}. \par

\vspace{3mm}

{\bf 8.}  Constants of Rosenthal - Dharmadhikari - Jogdeo - ...\par

 Let  $  p = \const \ge 2, \hspace{4mm}  \{ \zeta_k \} $ be a sequence of numerical centered, i.; i.d. r.v.  with finite $ p^{th} $ moment
 $ | \zeta|_p < \infty. $  The following constant,  more precisely, function on $ p, $ is called
 constants of Rosenthal-Dharmadhikari-Jogdeo-Johnson-Schechtman-Zinn-Latala-Ibragimov-Pinelis-Sharachmedov-Talagrand-Utev...:

$$
K_R(p) \stackrel{def}{=} \sup_{n \ge 1} \sup_{ \{\zeta_k\} } \left[ \frac{|n^{-1/2} \sum_{k=1}^n \zeta_k|_p}{|\zeta_1|_p} \right].  \eqno(1.10)
$$
 We will  use  the following ultimate up to an error value $ 0.5\cdot 10^{-5} $  estimate  for $ K_R(p), $ see \cite{Ostrovsky502} and reference therein:

 $$
 K_R(p) \le \frac{C_R \ p}{ e \cdot \log p}, \hspace{5mm}  C_R = \const := 1.77638.  \eqno(1.11)
 $$
 Note that for the symmetrical distributed r.v. $ \zeta_k $ the constant $  C_R $ may be reduced  up to a value $ 1.53572.$\par

\vspace{3mm}

{\bf 9.}  Estimates of maximum distribution for random fields,\\

\vspace{3mm}

 Let $ Y = Y( x,t), \ x \in X, \ t \in T  $ be separable  stochastic  continuous numerical  {\it random process (field)}
where in the capacity of "probability"  space is  $  (X, A, \mu); \ \overline{Y} := \sup_{t \in T} Y(t). $ \par

\vspace{3mm}

{\it  Notice that the measure $ \mu $ may be unbounded! } \par

\vspace{3mm}

 Let $  Q = \const \ge 1; $ we  denote

 $$
 \sigma(Q) = \sigma_Y(Q) = \sup_{t \in T} | Y(t, \cdot)|_Q = \sup_{t \in T}  \left[  \int_X |Y(t,x)|^Q \ \mu(dx)  \right]^{1/Q}
 \eqno(1.12)
 $$
 and suppose $ 0 <  \sigma(Q) < \infty. $
 Further, introduce the following natural  distance  $ d_Q  = d_Q(t,s), \ t,s \in T $ (more precisely,  semi - distance) on the set $ T $
as follows:

$$
 d_Q(t,s) \stackrel{def}{=} | Y(t,\cdot) - Y(s,\cdot)|_Q/\sigma(Q). \eqno(1.13)
$$

 Let $ \rho = \rho(t,s)  $ be arbitrary distance on the set $ T. $  We denote as usually by $ N(T,\rho, \epsilon) $ the
 minimal  number of closed $ \rho -  $ balls of radii $ \epsilon, \ \epsilon > 0 $ which cover the set $  T. $
 Evidently, $  \forall \epsilon  > 0 \ \Rightarrow N(T,\rho, \epsilon) < \infty $ iff  the set  $  T  $ is precompact set
 relative the distance $ \rho. $\par
  Recall that the quantity $   H(T,\rho, \epsilon) = \log  N(T,\rho, \epsilon)  $ is called "metric entropy of the set
$ T $  relative the distance $  \rho". $ \par
 For instance, if $ T $ is  closed bounded subset of the whole Euclidean space $  R^d $ containing a ball with positive radii,
 and

 $$
 \rho(t,s)  \asymp ||t-s||^{\alpha}, \ \alpha = \const \in (0, 1],
  $$
then $  N(T,\rho, \epsilon) \asymp \epsilon^{-d/\alpha}, \ \epsilon \in (0, \diam_{\rho}(T)/2) = (0, \radii(T)). $ \par
 Let $ \radii(T) = 1. $  It follows from theorem of Egorov that we can suppose the existence of the point $ t_0, \ t_0 \in T, $
 (non-random), "center of the set $ T", $ for which

 $$
 \sup_{t \in T} \rho(t, t_0) \le 1. \eqno(1.14)
 $$
 We agree $ N(T, \rho,1) = 1,  $  as long as the unit  ball with center in $  t_0 $  cover the set $  T. $ \par

\vspace{3mm}

  We refer  here the main results of articles  G.Pisier \cite{Pisier2} in the probabilistic case  $  \mu(X) = 1 $ and
 \cite{Ostrovsky303} more generally.

\vspace{3mm}

{\bf Proposition 1.1.} \\

\vspace{3mm}

$$
| \ \overline{Y} \ |_Q \le \sigma_Y(Q) \cdot \inf_{\theta \in (0,1)} \left[ \sum_{k=1}^{\infty} \theta^{k-1} \
N^{1/Q}(T, d_{Q}, \theta^k) \right]. \eqno(1.15)
$$

{\it Moreover, if the series in the  right - hand side (1.15) convergent, the r.f. $ Y = Y(t, \cdot) $
 is continuous almost everywhere relative the distance } $ d_Q: $

$$
\mu \{x: Y(\cdot,x) \notin  C(T, d_Q) \} = 0.
$$

\vspace{3mm}

 {\bf Example 1.1.}  Suppose

$$
N(T, d_Q, \epsilon) \le K^Q \ \epsilon^{-\kappa}, \ \epsilon \in (0,1), \ K = \const < \infty, \eqno(1.16)
$$
 where $ \kappa = \const \in [0,Q). $ The parameter $  \kappa, $ more precisely, its minimal  value,
 is called {\it entropic dimension } of the set $ T $ relative the distance $ d_Q. $\par

  We obtain after computations:

$$
| \ \overline{Y} \ |_Q \le  K \cdot \sigma_Y(Q) \cdot \left(1 - \frac{\kappa}{Q} \right)^{-1} \cdot \left[ \frac{\kappa}{Q} \right]^{-\kappa/(Q - \kappa)}.
\eqno(1.17)
$$

\vspace{3mm}

  Notice that the parameter $ \kappa $ may depend on the $ K,Q: \ \kappa = \kappa(K,Q). $ \par

\vspace{3mm}

{\bf Remark 1.1.}  Denote

$$
\nu(Q) = \sigma_Y(Q) \cdot \inf_{\theta \in (0,1)} \left[ \sum_{k=1}^{\infty} \theta^{k-1} \
N^{1/Q}(T, d_{Q}, \theta^k) \right], \eqno(1.18)
$$
and suppose that  there exists a value  $ Q_0 \in(1, \infty] $   for which   $ \forall Q < Q_0  \ \Rightarrow \nu(Q) < \infty.  $
 Further, denote
$$
h(Q) = Q \ \log \nu(Q), \ h^*(w) = \sup_{ Q \in (1, Q_0)}(w Q - h(Q)). \eqno(1.19)
$$
 It follows from Tchebychev's  inequality after optimisation over $  Q  $

$$
\mu \{x: (\overline{Y} > z) \} \le \exp \left( - h^*(\log z)  \right), \  z > 1. \eqno(1.20)
$$
 If in particular $ Q_0 = \infty, $  then from the estimate (1.20)  may be obtained the so - called {\it exponential} decreasing
bound  for tail of distribution of the value $ \overline{Y}, $ for example, of a view

$$
\mu \{x: (\overline{Y} > z) \} \le \exp \left( - C \ z^{\beta}  \right), \ \beta = \const > 0, \ z > 1. \eqno(1.21)
$$
 See for detail explanation \cite{Ostrovsky1}, chapter 1, sections 1.1 - 1.5. \par

\vspace{3mm}

\section{ Moment estimates in the first norm.}

\vspace{3mm}

{\it  In  what follows } $ \xi = \xi(x,t) = \xi(\omega; x,t), \ \omega \in \Omega, x \in X, \ t \in T $  {\it be measurable
separable stochastic continuous  numerical random field (r.f.),   }  $ p \ge 2, \ \xi_i = \xi_i(x,t) = \xi_i(\omega; x,t), \ i=1,2,\ldots $
{\it  be independent copies of    } $ \xi(x,t), $\par

$$
|\xi|_{p,\infty} = |\xi|_{p,X;\infty,T}  = \sup_{t \in T}  \left[ \int_X |\xi(x,t)|^p \ \mu(dx) \right]^{1/p}. \eqno(2.1)
$$

\vspace{3mm}

 We define an auxiliary   random field

 $$
 \eta(t)  = \eta_p(t) = \int_X |\xi(x,t)|^p \ \mu(dx), \eqno(2.2)
 $$
so that

$$
|\xi|^p_{p,\infty} = \sup_{t \in T} \eta_p(t). \eqno(2.3)
$$

We intend ro use the Pisier's  estimate (1.15) or its generalization for  the non-norming measure \cite{Ostrovsky303}.
For this purpose we need to make some calculations. \\
 First of all we need to estimate the value

 $$
 \sigma_{p,Q} := \sup_{t \in T} | \eta_p(t)|_{Q, \Omega}.  \eqno(2.4)
 $$
 We have using Minkowsky and permutation inequalities:

 $$
 | \eta_p(t)|_{Q, \Omega} \le \int_X \left| \ | \xi(x,t) |^p \   \right|_{Q, \Omega} \ \mu(dx) =
 \int_X \left[ {\bf E} |\xi(x,t)|^{ p Q} \right]^{1/Q} \ \mu(dx),
 $$
therefore

$$
\sigma_{p,Q} \le \overline{\sigma}_{p,Q} \stackrel{def}{=} \sup_{t \in T} \int_X \left[ {\bf E} |\xi(x,t)|^{p Q} \right]^{1/Q}  \ \mu(dx). \eqno(2.5)
$$

 Further, we need to estimate the $  L_Q(\Omega) $ norm of an increments $ \eta_p(t) - \eta_p(s). $ Since

$$
\eta_p(t) - \eta_p(s) = \int_X \left[ |\xi(x,t)|^p - |\xi(s,x)|^p      \right] \ \mu(dx),
$$
we  deduce analogously to the inequality (2.5):

$$
d_{p,Q}(t,s) \stackrel{def}{=} |\eta_p(t) - \eta_p(s)|_{Q, \Omega} \le
$$

$$
\int_X \left| {\bf E} |\xi(x,t)|^p - {\bf E}|\xi(x,s)|^p  \right|^{1/Q} \mu(dx) =: \overline{d}_{p,Q}(t,s). \eqno(2.6)
$$

 If we define the norming distance

 $$
 \rho_{p,Q}(t,s) := \frac{d_{p,Q}(t,s)}{ \overline{\sigma}_{p,Q}}, \eqno(2.7)
 $$
then

$$
N(T, \rho, \epsilon) = N(T, d_{p,Q}, \epsilon \ \overline{\sigma}_{p,Q}) \le  N(T, \overline{d}_{p,Q}, \epsilon \ \overline{\sigma}_{p,Q}).
$$
 We obtain using the inequality (1.15):

 \vspace{3mm}

 {\bf Proposition 2.1.}   {\it Let $ p \ge 2, \ Q \ge 1. $ Define the following function:}

 $$
 \psi_p^p(Q) \stackrel{def}{=}   \overline{\sigma}_{p,Q} \inf_{\theta \in (0,1)}
\left[  \sum_{k=1}^{\infty} \theta^{k-1} N^{1/Q} \left(T, \overline{d}_{p,Q}, (\theta \ \overline{\sigma}_{p,Q})^k \right) \right]. \eqno(2.8)
 $$
{\it  and suppose it finiteness for some values $ Q; $ otherwise it is nothing to prove. Assertion: }

$$
\left\{ {\bf E} |\xi(\cdot, \cdot)|^{pQ}_{p,\infty} \right\}^{1/p Q} \le \psi_p(Q). \eqno(2.9)
$$

\vspace{3mm}

{\bf Example 2.1.} Suppose for some positive finite constants $ c_1(p), \ m(p) $ and for all the values   $  Q, \ Q \ge 1 $

$$
\psi_p^p(Q) \le c_1(p) \ Q^{m(p)}.
$$
 We deduce from (2.9) the following {\it exponential }  tail estimate:

 $$
{\bf P} \left( |\xi(\cdot, \cdot)|^{p}_{p,\infty} > x \right)  \le \exp \left( - c_2(p) x^{1/m(p)} \right), \ x > 0. \eqno(2.10)
 $$

\vspace{3mm}

\section{ Central Limit Theorem and tail estimates \\
 for normed sums of random vectors\\
 in the first norm.}

\vspace{3mm}

 {\it Assume in addition to the second section that  $  p \ge 2 $ and that the r.f. } $  \xi(x,t),  $ {\it and with it the r.f. }
$  \xi_i(x,t)  $ {\it  are mean zero: } $ {\bf E} \xi_i(x,t) = 0, $ and denote

$$
S_n(x,t) = n^{-1/2} \sum_{i=1}^n \xi_i(x,t). \eqno(3.1)
$$

$$
\tau_p^{(n)}(t) = \int_X | S_n(x,t)|^p \ \mu(dx) = |S_n(\cdot,t)|^p_{p,X}. \eqno(3.2)
$$

\vspace{3mm}

{\bf 1.}  We intend to  estimate  uniformly over numbers of summand $ n $  first of all in this section
the moments of the random variable

$$
\zeta_p = \sup_{t \in T} \tau_p^{(n)}(t), \eqno(3.3)
$$
i.e.  the values

$$
g_p(Q) = \sup_n \left| \ \sup_{t \in T}  \tau_p^{(n)}(t) \  \right|_{Q, \Omega}.
$$

 Note that

 $$
 \left[ g_p(Q) \right]^{1/p} =  \sup_n | S_n(\cdot,\cdot)|_{p,X; \infty,T; Q,\Omega}.
 $$

\vspace{3mm}

 Some new notations:  $ \rho_{v,x}(t,s) : = $

 $$
 | \ \xi(x,t) - \xi(x,s) \ |_{v,\Omega} =  \left\{ \left[  {\bf E} |\xi(x,t) - \xi(x,s)| \right]^v    \right\}^{1/v},
\ v = \const \ge 1; \eqno(3.4),
 $$

$$
W_{\gamma}(x) = \sup_{t \in T} | \ \xi(x,t) \ |_{\gamma, \Omega} = \sup_{t \in T} \left\{ {\bf E} |\xi(x,t)|^{\gamma} \right\}^{1/\gamma}, \eqno(3.5)
$$

$$
J(t,s; p,Q; \alpha,\beta) = \int_X W^{p-1}_{ (p-1) \beta Q }(x) \ \rho_{ \alpha Q,x }(t,s) \ \mu(dx), \ \alpha,\beta > 1,
1/\alpha + 1/\beta = 1;
$$

$$
r_{p,Q}(t,s) = 2 \ p \ \inf_{\alpha, \beta} \left[ K_R(\alpha Q) \ K_R^{p-1}((p-1) \beta Q) \ J(t,s; p,Q; \alpha,\beta)  \right]. \eqno(3.6)
$$

 Evidently, $ r_{p,Q}(t,s)  $ is the distance as the function on $ (t,s),  $ if it is finite.
 The minimum in the right - hand side (3.6) is calculated over all the values $ (\alpha, \beta) $ for which
 $  \alpha, \beta > 1, \ 1/\alpha + 1/\beta = 1.  $ \par

 Further, denote

 $$
 \hat{\sigma}_{p,Q}:= K_R^p(p \ Q) \ \overline{\sigma}_{p,Q}, \eqno(3.7)
 $$

$$
\hat{r}_{p,Q}(t,s) := r_{p,Q}(t,s)/\hat{\sigma}_{p,Q},\eqno(3.8)
$$

$$
 \nu_p^p(Q) \stackrel{def}{=}   \hat{\sigma}_{p,Q} \cdot \inf_{\theta \in (0,1)}
\left[  \sum_{k=1}^{\infty} \theta^{k-1} N^{1/Q} \left(T, \hat{r}_{p,Q}, (\theta \ \hat{\sigma}_{p,Q})^k \right) \right]. \eqno(3.9)
 $$

\vspace{4mm}

{\bf Theorem 3.1.} {\it If }   $  \nu_p(1) < \infty,  $  {\it then    } $  \{ \xi_i(x,t)  \} $ {\it  satisfies the CLT in the space
 $ CL_p(X,T). $  }\\

\vspace{4mm}

{\bf Theorem 3.2.} {\it If for some}   $ Q = \const \ge 1 \ \Rightarrow  \nu_p(Q) < \infty,  $  {\it then    }

\vspace{3mm}

$$
 \sup_n \left\{ {\bf E} |S_n(\cdot, \cdot)|^{pQ}_{p,\infty} \right\}^{1/p Q} \le \nu_p(Q). \eqno(3.10)
$$

\vspace{4mm}

{\bf Proofs.}\\
{\bf 1.} We need first of all to obtain the estimate (3.10), of course,  through the proposition 2.1. We have
using the Rosenthal's constants and the Minkowsky  inequality:

$$
| \ \tau_p^{(n)}(t) \ |_{Q, \Omega}  = \left| \ \int_X |S_n(x,t)|^p \ \mu(dx) \ \right|_{Q, \Omega} \le
$$

$$
 \int_X  \left| \ |S_n(x,t)|^p  \ \right|_{Q, \Omega}   \ \mu(dx)  \le
 \int_X K_R^p(p \ Q) \ \left| \ |\xi(x,t)|^p  \ \right| _{Q, \Omega}  \ \mu(dx)\le
$$

$$
K_R^p(p \ Q) \ \overline{\sigma}_{p, Q} =  \hat{\sigma}_{p, Q}. \eqno(3.11)
$$

\vspace{3mm}
{\bf 2.} The estimation of a difference
$$
 \Delta \tau(t,s) =  \tau_p^{(n)}(t) - \tau_p^{(n)}(s)
$$
is more complicated. We have consequently:

$$
\Delta \tau = \int_X \left[ |S_n(x,t)|^p - |S_n(x,s)|^p  \right]  \ \mu(dx),
$$

$$
|\Delta \tau |_{Q,\Omega} \le \int_X \left| \ |S_n(x,t)|^p - |S_n(x,s)|^p \  \right|_{Q, \Omega}  \ \mu(dx) =
$$

$$
\int_X \left[ {\bf E} \left| \ |S_n(x,t)|^p -  |S_n(x,s)|^p  \ \right|^Q  \right]^{1/Q} \ \mu(dx). \eqno(3.12)
$$

  We exploit the following  elementary inequality:

$$
| \ |x|^p - |y|^p \ | \le p \cdot |x-y| \cdot \left[ |x|^{p-1} + |y|^{p-1} \right], \ x,y \in R, \eqno(3.13)
$$
and obtain after substituting into (3.12), where $  x = S_n(x,t), \ y = S_n(x,s):  \ |\Delta \tau |_{Q,\Omega}/p \le  $

$$
 \int_X \left| \ |S_n(x,t) - S_n(x,s)| \ \cdot
 \ \left[ \ |S_n(x,t)|^{p-1} + |S_n(x,s)|^{p-1} \right] \ \right|_{Q, \Omega}  \ \mu(dx). \eqno(3.14)
$$

  It follows from the H\"older's inequality

$$
|\eta_1 \eta_2|_{Q,\Omega} \le |\eta_1|_{\alpha Q, \Omega} \cdot |\eta_2|_{\beta Q, \Omega},
$$
where as before $ \alpha, \beta > 1, \ 1/\alpha + 1/\beta = 1. $ Therefore

$$
|\Delta \tau |_{Q,\Omega}/p \le \int_X g_1(t,s; \alpha, x) \cdot g_2(\beta,x) \ \mu(dx),
$$
where

$$
g_1(t,s; \alpha, x) =  g_1(t,s; \alpha, x; \Omega)  = |  S_n(x,t) - S_n(x,s) |_{\alpha Q; \Omega}, \eqno(3.15)
$$
and

$$
g_2(\beta,x) = g_2(\beta,x; p,\Omega) = \sup_{t,s \in T} \left| \ \left[|S_n(x,t)|^{p-1} + |S_n(x,s)|^{p-1} \right] \ \right|_{\beta Q, \Omega} . \eqno(3.16)
$$

 We estimate $ g_1(\cdot) $ using Rosenthal's inequality:

 $$
 g_1(t,s; \alpha, x; \Omega) \le K_R(\alpha Q) \ | \xi(x,s) - \xi(x,s) |_{\alpha Q}= \rho_{\alpha Q,x}(t,s). \eqno(3.17)
 $$
 Further,

 $$
 g_2(\beta,x; p,\Omega)  \le 2 \ K_R^{p-1}(\beta (p-1))  \sup_{t \in T}  | \xi(x,t)|^{p-1}_{\beta(p-1), \Omega} =
 $$

$$
 2 \ K_R^{p-1}(\beta (p-1))\ W_{\beta(p-1)}^{p-1}(x). \eqno(3.17)
$$
 We get  after substituting into (3.15) and (3.16)

$$
|\Delta \tau |_{Q,\Omega} \le r_{p,Q}(t,s). \eqno(3.18)
$$

 It remains to use proposition 2.1. Theorem 3.2 is proved. \\

\vspace{3mm}

{\bf 3. Proof of theorem 3.1} is similar to one for the mixed spaces in \cite{Ostrovsky303}.     \par

 Let $  \gamma_p(1) < \infty. $ By theorem 3.2

$$
 \sup_n \left\{ {\bf E} |S_n(\cdot, \cdot)|^{p}_{p,\infty} \right\}^{1/p} \le \nu_p(1) < \infty. \eqno(3.19)
$$

  As long as the Banach space $ CL_p(X) $  is separable
and the function $  y \to |y|^{p}  $ satisfies the $ \Delta_2  $ condition, there exists a linear compact operator
$  U: \ CL_p(X) \to CL_p(X)  $ such that

$$
{\bf P} \left(U^{-1} S_n(\cdot, \cdot) \in CL_p(X) \right) = 1 \eqno(3.20)
$$
and moreover

$$
 \sup_n  {\bf E} | U^{-1} S_n |^p_{p,\infty} < \infty. \eqno(3.21)
$$
 \cite{Ostrovsky2}; see also \cite{Buldygin1}, \cite{Ostrovsky603}.\par

 We get using Tchebychev's  inequality

 $$
\sup_n {\bf P} \left( |U^{-1}[ S_n]|_{p,X} > Z   \right) \le C(p)/Z^p < \epsilon,\eqno(3.22)
 $$
for sufficiently greatest values $ Z = Z(\epsilon), \ \epsilon \in (0,1). $ \par
  Denote by $  W  = W(Z) $ the set

$$
W = \{ f: f \in L_p(X),  |U^{-1}[f]|_{p,X}  \le Z \}. \eqno(3.22)
$$
 Since the operator $ U $ is compact, the set $  W  = W(Z) $ is compact set in the space $  CL_p(X). $ It follows from inequality
(3.22) that

$$
\sup_n {\bf P} \left( S(n) \notin W(Z) \right) \le \epsilon.
$$
  Thus, the sequence $ \{  S_n \}  $ satisfies the famous Prokhorov's criterion \cite{Prokhorov1} for  weak compactness of
the family of distributions in the separable metric spaces. \par
 This completes the proof of theorem 3.1. \par

\vspace{3mm}

{\bf Remark 3.1.}  Another way to prove the   the weak compactness of distributions $  S_n(\cdot, \cdot) $ is
 by using theorem 4.3.2  from monograph \cite{Ostrovsky1}, chapter 4, section 3. \\
 Actually, if $ \nu_p(1) < \infty, $ then $ \forall h > 0 \ \Rightarrow $

 $$
 \lim_{\epsilon \to 0+} \sup_n
 {\bf P} \left( \sup_{\rho(t,s) < \epsilon} \left|\tau_p^{(n)}(t) - \tau_p^{(n)}(s) \right| > h \right) = 0. \eqno(3.23)
 $$
 The proposition of theorem 3.1 follows from  theorem 1 of monograph \cite{Gikhman1}, p.408.\\

\vspace{3mm}

\section{ Moment estimations in the second norm.}

\vspace{3mm}

 Recall that for the  measurable numerical function $  f = f(x,t), \ x \in X, t \in T $

 $$
 | f(\cdot, \cdot) |_{\infty,p} = | f(\cdot, \cdot) |_{\infty,T;p,X} = \left[ \int_X \sup_{t \in T} |f(x,t)|^p \ \mu(dx)   \right]^{1/p},
 \eqno(4.0)
 $$
and for the function $ f \in L(p)C $

$$
\lim_{\epsilon \to 0+} \left[ \int_X \sup_{d(t,s) < \epsilon } |f(x,t) - f(x,s)|^p \ \mu(dx) \right]^{1/p} = 0.
$$

\vspace{4mm}

{\bf 1.} Let $ \xi = \xi(x,t) = \xi(\omega; x,t)  $ be triple measurable numerical random field; then

 $$
 |\xi|_{\infty,p}^p = \int_X \sup_{t \in T} |\xi(x,t)|^p \ \mu(dx),
 $$
and if we denote  again

$$
\eta_p(t) = |\xi(x,t)|^p, \hspace{5mm} \Delta_{p,Q} = \sup_{t \in T} | \eta_p(t)|_{Q,X};
$$

$$
\Delta_{p,Q} =  \sup_{t \in T} \left[ \int_X |\xi(x,t)|^{p Q} \ \mu(dx)   \right]^{1/Q}.  \eqno(4.1)
$$

 Note that $ \Delta_{p,Q} $ is a random variable.\par

\vspace{3mm}


{\bf 2.}  We introduce a {\it random} distance on the set $  T  $ as follows:

$$
\delta_{p,Q}(t,s) = | \ \eta_p(t) - \eta_p(s) \ |_{Q,X}/\Delta_{p,Q} =
$$

$$
\sqrt[Q]{\int_X \left| \ |\xi(x,t)|^p - |\xi(x,s)|^p  \right|^Q \ \mu(dx)}/\Delta_{p,Q}, \eqno(4.2)
$$
and the following random entropy function:

$$
\lambda_p(Q) = \Delta_{p,Q} \cdot \inf_{\theta \in (0,1)} \left[ \sum_{k=1}^{\infty} \theta^{k-1} \ N(T, \delta_{p,Q}, \theta^k)   \right].
\eqno(4.3)
$$

\vspace{3mm}

{\bf 3. Proposition 4.1.}

$$
\sqrt[pQ]{ {\bf E} |\xi|^{pQ}_{p,X; \infty,T}} \ = \ |\xi|_{p,X; \infty,T; pQ,\Omega} \ \le
\ \sqrt[p Q]{ {\bf E} \lambda_p^Q(Q)}. \eqno(4.4)
$$

\vspace{3mm}

{\bf Proof } follows immediately from the proposition 2.1,  in which in the capacity the probability space  used the
measurable space $ (X, A, \mu)  $ with non - normed measure $ \mu. $ \par

Of course, this statement (4.4) is meaningful  only for those values $  Q, $ for which
 $ {\bf E} \lambda_p^Q(Q) < \infty. $ \par

\vspace{3mm}

\section{ Concluding remarks.}

\vspace{3mm}

{\bf  CLT for dependent r.v. in hybrid spaces.}\\

\vspace{3mm}
 Analogously to the article \cite{Ostrovsky303} may be considered the case when the r.v. $ \xi_i = \xi_i(x,t) $ dependent,
for example, form a {\bf  martingale } or {\bf mixingale} sequence. \\

\vspace{3mm}

{\bf  Martingale case.} \\

\vspace{3mm}

 We suppose for example as before that  $ \{ \xi_k(\cdot)  \}  $ are mean zero and form a strictly stationary sequence,
 $ p \ge 2. $
 Assume in addition that $ \{ \xi_k(\cdot)  \}  $ form a martingale difference sequence
 relative certain filtration $ \{  F(k) \},  \ F(0) = \{ \emptyset, \Omega\},  $
 $$
  {\bf E} \xi_k/F(k) = \xi_k, \hspace{5mm}  {\bf E} \xi_k/F(k-1) = 0, \ k= 1,2,\ldots.
 $$

 Then the proposition  of theorem 3.1, 3.2 remains true; the estimate of theorem 3.2 is also true  up to multiplicative
absolute constant. \par
 Actually, the convergence of correspondent characteristical functionals follows from the ordinary one - dimensional
 CLT for martingales, see in the classical monograph of  Hall P., Heyde C.C.
  \cite{Hall1}, chapter 2; the Rosenthal's constant for the sums  of martingale differences with at the same up to
  multiplicative constant coefficient is obtained by A.Osekowski \cite{Osekowski1}, \cite{Osekowski2}.
 See also \cite{Ostrovsky3}.\par

\vspace{3mm}

{\bf  Mixingale case.} \\

\vspace{3mm}

 We suppose again  that  $ \{ \xi_k(\cdot)  \}  $ are mean zero and form a strictly stationary sequence,
 $  p \ge 2. $  This  sequence is said to be {\it mixingale, } in the terminology of the book
\cite{Hall1}, if it satisfies this or that mixing condition.\par

 We  consider here only the  superstrong mixingale. Recall that the superstrong, or $  \beta = \beta(F_1, F_2) $
index between two sigma-algebras  is defined as follows:

$$
\beta(F_1, F_2) = \sup_{A \in F_1, B \in F_2, {\bf P}(A) {\bf P}(B) > 0 } \left| \frac{{\bf P}(AB) - {\bf P}(A) {\bf P}(B)}{{\bf P}(A) {\bf P}(B)} \right|.
$$

 Denote

 $$
 F_{-\infty}^0  = \sigma(\xi_s, \ s \le 0),  \hspace{5mm} F_n^{\infty} = \sigma(\xi_s, \ s \ge n),
 $$

$$
\beta(n) = \beta \left(F_{-\infty}^0 , F_n^{\infty} \right),
$$

 The  sequence $ \{\xi_k \} $ is said to be {\it superstrong mixingale, } if $ \lim_{n \to \infty} \beta(n) = 0. $ \par
This notion  with some applications was introduced and investigated by B.S.Nachapetyan and R.Filips \cite{Nachapetyan1}.
See also   \cite{Ostrovsky3}, \cite{Ostrovsky1}, p. 84 - 90. \par

\vspace{3mm}

 Introduce the so-called mixingale Rosenthal coefficient:

 $$
 K_M(m) =  m \ \left[  \sum_{k=1}^{\infty} \beta(k) \ (k+1)^{ (m -2)/2  }   \right]^{1/m}, \ m \ge 1.
 $$

B.S.Nachapetyan in \cite{Nachapetyan1} proved that for the superstrong centered  strong stationary strong mixingale
sequence $  \{ \eta_k \} $ with $ K_M(m) < \infty $ the following estimate is true:

$$
\sup_{n \ge 1}  \left| n^{-1/2} \sum_{k=1}^n \eta_k \right|_m \le C \cdot K_M(m) \cdot |\eta_1|_m,
$$
so that the "constant" $ K_M(m) $ play at the same role for mixingale as the Rosenthal  constant  for independent variables.\par

 As a consequence:  theorems 3.1 and 3.2  remains true for strong mixingale sequence $ \{  \xi_k \}: $
theorem 3.1 under conditions: $ K_M (p) < \infty $  for theorem 3.1 and $ K_M(pQ) < \infty $  for the theorem 3.2
with replacing $  K_R(pQ) $ on the expression $  K_M( pQ ). $  \par

\vspace{3mm}

{\bf  Another approach.}

\vspace{3mm}

Another approach for the tail estimation for the maximum distribution of random field is closely related with notions
"majorizing measures" or equally "generic chaining", see \cite{Fernique1}, \cite{Talagrand1},  \cite{Talagrand2},  \cite{Talagrand3},
\cite{Bednorz1},  \cite{Bednorz2}, \cite{Ostrovsky404}, \cite{Ostrovsky405} etc. \par
   But by our opinion offered here method is more convenient for the announcement goals. \par


\vspace{6mm}
\begin{thebibliography}{99}

\vspace{4mm}

\bibitem{Adams1}
{\sc Adams R.A.} {\it  Anisotropic Sobolev Inequalities.}
Casopic pro Pestovani Matematiky, (Prague),  No. 3, 267—279.

\bibitem{Bednorz1}
{\sc Bednorz W.} (2006). {\it A theorem on Majorizing Measures.} Ann. Probab., 34,
1771-1781. MR1825156.

\bibitem{Bednorz2}
{\sc Bednorz W.} {\it The majorizing measure approach to the sample boundedness.}
arXiv:1211.3898v1 [math.PR] 16 Nov 2012

\bibitem{Benedek1}
{\sc Benedek A. and Panzone  R.} {\it The space $  L_p $ with mixed norm.} Duke Math. J., {\bf 28}, (1961),  301 - 324.

\bibitem{Besov1}
{\sc Besov O.V., Il’in V.P., Nikol’skii S.M. } {\it Integral representation of functions
and imbedding theorems.} Vol.1; Scripta Series in Math., V.H.Winston
and Sons, (1979), New York, Toronto, Ontario, London.

\bibitem{Buldygin1}
{\sc Buldygin V.V.} (1984). {\it Supports of probabilistic measures in separable Banach
spaces.} Theory Probab. Appl. 29 v.3, pp. 528 - 532, (in Russian).

\bibitem{Hall1}
{\sc Hall P., Heyde C.C. } {\it Martingale Limit Theory and Applications.}  Academic
Press, New York. (1980)

\bibitem{Fernique1}
{\sc Fernique X.} {\it Regularite de fonctions aleatoires non gaussiennes.} Ecolee de
Ete de Probabilits de Saint-Flour XI-1981. Lecture Notes in Math. 976, (1983),
174, Springer, Berlin.

\bibitem{Fujita1}
{\sc  Fujita H. and Kato T. } {\it On the Navier-Stokes initial value problem I. } Arch. Ration. Mech.
Anal., 16(1964), 269 \ – \ 315.

\bibitem{Gikhman1}
{\sc Gikhman I.I., Skorokhod A.V.} {\it The Theory of Stochastic Processes, I; } Springer Verlag, (1980),
Berlin - Heidelberg - New York.

\bibitem{Kato1}
{\sc  Kato T. }  {\it Strong $ L_p $ solutions of the Navier-Stokes equations in $ R^m $ with applications to
weak solutions. } Math. Zeitschrift, {\bf 187,} \ (1984), 471 \ – \ 480.

\bibitem{Kato2}
{\sc Kato T. and  Ponce G.} {\it Commutator estimates and the Euler and Navier-Stokes equations. }
Comm. P. D. E., 41(1988), 891 \ - \ 907.

\bibitem{Kufner1}
{\sc  Kufner A., Oldrich J., and Fucik S.} {\it  Function Spaces.}  Noordhoff International  Publishing,(1977),
Academia Publishing House  of the Chechoslovak Academy of Science,  Prague.

\bibitem{Lions1}
{\sc Lions J.L.} {\it Quelques methodes de  resolutions des problemes auxlimites non linearies. } Dunod;
Gauthier - Villars, Paris, (1969), MR 41.

\bibitem{Lions2}
{\sc Lions J.L., Magenes E.}  {\it  Non - homogeneous boundary value problems and applications.  } (1972)
V.1, Springer Verlag,  Berlin - Heidelberg - New York, MR 40.

\bibitem{Nachapetyan1}
{\sc Nachapetyan B.S.} {\it On the certain criterion of weak dependence. } Probab. Theory Appl., (1980),
{\bf 2, } V. 26, 374 - 381.

\bibitem{Leoni1}
{\sc Leoni G. } {\it A first Course in Sobolev Spaces.} Graduate Studies in Mathematics,
v. 105, AMS, Providence, Rhode Island, (2009).

\bibitem{Lieb1}
{\sc Lieb E., Loss M.} {\it Analysis.} Providence, Rhode Island, 1997.

\bibitem{Ostrovsky1}
{\sc  Ostrovsky E.I.} (1999). {\it Exponential estimations for random Fields and its
Applications, (in Russian).}  Moscow - Obninsk, OINPE.

\bibitem{Ostrovsky2}
{\sc Ostrovsky E.I.} (1980).{\it On the support of probabilistic measures in separable
Banach spaces.} Soviet Mathematic, Doklady, v.255, No 6 pp. 836 - 838, (in Russian).

\bibitem{Osekowski1}
{\sc Osekowski A.} {\it Inequalities for dominated martingales.} Bernoulli 13 (2007), 54-79.

\bibitem{Osekowski2}
{\sc Osekowski A.}  {\it Sharp martingale and semimartingale inequalities.} Monografie Matematyczne 72,
Birkhauser, 2012.

\bibitem{Ostrovsky3}
{\sc  Ostrovsky E. and Sirota L.} {\it  Moment and tail estimates for martingales and martingale transform,
with application to the martingale limit theorem in Banach spaces.}
arXiv:1206.4964v1 [math.PR] 21 Jun 2012

\bibitem{Ostrovsky603}
{\sc Ostrovsky E.} {\it Support of Borelian measures in separable Banach spaces. }
arXiv:0808.3248v1 [math.FA] 24 Aug 2008

\bibitem{Ostrovsky502}
{\sc  Ostrovsky E. and Sirota L.}  {\it  Schl\"omilch and Bell series for Bessel's functions, with
probabilistic applications.}
arXiv:0804.0089v1 [math.CV] 1 Apr 2008

\bibitem{Talagrand1}
 {\sc Talagrand M.} (1996). {\it Majorizing measure: The generic chaining.}
 Ann. Probab., {\bf 24} 1049 - 1103. MR1825156

\bibitem{Talagrand2}
 {\sc Talagrand M.} (2005). {\it The Generic Chaining. Upper and
     Lower Bounds of Stochastic Processes.} Springer, Berlin. MR2133757.

\bibitem{Talagrand3}
{\sc Talagrand M.} (1990). {\it Sample boundedness of stochastic processes under
increment conditions. } Ann. Probab. 18, 1 - 49.

\bibitem{Taylor1}
{\sc Taylor M.E.} {\it Partial Differential Equations III. Non-linear Rquations.} Applied
Math. Sciencies, 117, Springer, (1996).

\vspace{16mm}

\bibitem{Araujo1}
{\sc Araujo  A., Gine E. }
{\it The central limit theorem for real and Banach valued random variables.}
Wiley, (1980), London, New York.

\bibitem{Billingsley1}
{\sc Billingsley P.} {\it Probability and measure.}
Wiley, 1979, London, New York.

\bibitem{Billingsley2}
{\sc Billingsley P.} {\it  Convergence of probability measures.}
Wiley, (1968), London, New York.

\bibitem{Dudley1}
{\sc Dudley R.M.} {\it Uniform Central Limit Theorem}. Cambridge University Press, (1999)

\bibitem{Grenander1}
{\sc Grenander U.} {\it Probabilities on algebraic structures.}
Wiley, 1963; London, New York.

\bibitem{Ledoux1}
 {\sc Ledoux M., Talagrand M.} (1991) {\it Probability in Banach Spaces.}
      Springer, Berlin, MR 1102015.

\bibitem{Ostrovsky1}
{\sc  Ostrovsky E.I.} (1999). {\it Exponential estimations for random Fields and its
applications (in Russian).}  Moscow - Obninsk, OINPE.

\vspace{16mm}

\bibitem{Fortet1}
{\sc Fortet R. and Mourier E.} {\it Les  fonctions  alratoires comme elements aleatoires dans
les espaces de Banach.} Studia Math., {\bf 15}, (1955), 62-79.

\bibitem{Garling1}
{\sc Garling D.J.H. }
{\it Functional Central Limit Theorems in Banach Spaces.}
 The Annals of Probability, Vol. 4, No. 4 (Aug., 1976), pp. 600-611

\bibitem{Gine1}
{\sc Gine E.} {\it On the Central Limit theorem for sample continuous processes.} Ann.
Probab. (1974), 2, 629-641.

\bibitem{Gine2}
{\sc Gine E., Zinn J.} {\it  Central Limit Theorem and Weak Laws of Large Numbers in certain Banach Spaces.  }
Z. Wahrscheinlichkeitstheory verw. Gebiete. {\bf 62}, (1983), 323  -  354.

\bibitem{Heinkel1}
{\sc Heinkel B.} Measures majorantes et le theoreme de la limite centrale dans
C(S). Z. Wahrscheinlichkeitstheory. verw. Geb., (1977). 38, 339-351.

\bibitem{Jain1}
{\sc Jain N.C. and Marcus M.B.} {\it Central limit theorem for $C(S)$ valued random
variables.} J. of Funct. Anal., (1975), 19, 216-231.

\bibitem{Kozachenko1}
 {\sc Kozachenko Yu. V., Ostrovsky E.I.} (1985). {\it The Banach Spaces of
      random Variables of subgaussian type.} Theory of Probab. and Math.
      Stat. (in Russian). Kiev, KSU, {\bf 32}, 43 - 57.

\bibitem{Ostrovsky301}
{\sc Ostrovsky E., L.Sirota L.}
{\it CLT for continuous random processes under approximations terms.}
arXiv:1304.0250v1 [math.PR] 31 Mar 2013


\bibitem{Ostrovsky302}
{\sc Ostrovsky E., Rogover E.} {\it  Maximal inequalities in bilateral Grand Lebesgue Spaces.  }
arXiv:0808.3247v1 [math.FA] 24 Aug 2008

\bibitem{Ostrovsky303}
{\sc Ostrovsky E., L.Sirota L.} {\it Central Limit Theorem and exponential tail estimates in mixed
(anosotropic) Lebesgue spaces.  }
arXiv:1308.5606v1 [math.PR] 26 Aug 2013

\bibitem{Ostrovsky304}
{\sc Ostrovsky E.} {\it Exact exponential estimations for random field maximum
distribution.} (2002), Theory Probab. Appl. 45 v.3, 281 - 286.


\bibitem{Pisier1}
{\sc  Pisier G.,  Zinn J.} {\it On the limit theorems for random variables with values in the spaces } $ L_p, \ 2  \le p  < \infty. $
 Z. Wahrscheinlichkeitstheorie verw. Gebiete 41, 289 - 304 (1978).

\bibitem{Pisier2}
{\sc Pizier G.} {\it Condition d'entropic assupant la continuite de certain processus
et applications a l'analyse harmonique.} Seminaire d'analyse fonctionnalle. (1980)
Exp. 13 p. 23 - 24.

\bibitem{Prokhorov1}
{\sc Prokhorov Yu.V.} {\it Convergense of Random Processes and Limit Theorems
of Probability Theory.} Probab. Theory Appl., (1956), V. 1, 177-238.

\bibitem{Rackauskas1}
{\sc Rackauskas A, Suquet Ch.}
{\it Central limit theorems in H\"ölder topologies for Banach space valued random fields.}
Teor. Veroyatnost. i Primenen., 2004, Volume 49, Issue 1, Pages 109–125 (Mi tvp238)

\bibitem{Rhee1}
{\sc Rhee Wan Soo and Michel Talagrand M. }
{\it Uniform bound in the central limit theorem for Banach space valued dependent random variables},
 Journal of Multivariate Analysis, 1986, vol. 20, issue 2, pages 303-320

\bibitem{Song1}
{\sc Song L.} {\it A counterexample in the Central Limit Theorem}.
Bulletin of the London Mathematical Society, Volume 31,
 Issue 02, March 1999, pp 222-230

\bibitem{Sualb1}
{\sc Sualb Z.} {\it Central limit theorems for random processes with
sample paths in exponential Orlicz spaces.}
 Stochastic Processes and their Applications 66, (1997), l-20.

\bibitem{Zinn1}
{\sc Zinn J.  } {\it A Note on the Central Limit Theorem in Banach Spaces.}
 Ann. Probab. Volume 5, Number 2 (1977), 283-286.

\vspace{16mm}

\bibitem{Frolov1}
{\sc Frolov A.S., Tchentzov N.N. } {\it On the calculation by the Monte-Carlo
method definite integrals depending on the parameters. } Journal of Computational
Mathematics and Mathematical Physics, (1962), V. 2, Issue 4, p. 714-718 (in
Russian).

\bibitem{Grigorjeva1}
{\sc Grigorjeva M.L., Ostrovsky E.I.} {\it Calculation of Integrals on discontinuous
Functions by means of depending trials method.} Journal of Computational
Mathematics and Mathematical Physics, (1996), V. 36, Issue 12, p. 28-39 (in
Russian).

\bibitem{Ostrovsky402}
{\sc Ostrovsky E., Sirota L.} {\it Monte-Carlo method for multiple parametric integrals
calculation and solving of linear integral Fredholm equations of a second
kind, with confidence regions in uniform norm.}
 arXiv:1101.5381v1 [math.FA] 27 Jan 2011

\bibitem{Ostrovsky403}
{\sc  Ostrovsky E., Rogover E.} {\it Non - asymptotic exponential bounds for
MLE deviation under minimal conditions via classical and generic chaining methods.}
arXiv:0903.4062v1 [math.PR] 24 Mar 2009

\bibitem{Ostrovsky404}
{\sc Ostrovsky E., Sirota L.} {\it  Simplification of the majorizing measure method,
with development.}
arXiv:1302.3202v1 [math.PR] 13 Feb 2013

\bibitem{Ostrovsky405}
{\sc  Ostrovsky E., Rogover E.}
{\it Exact exponential Bounds for the random Field
Maximum Distribution via the Majorizing Measures (Generic Chaining.)}
 arXiv:0802.0349v1 [math.PR] 4 Feb 2008

\vspace{4mm}

\end{thebibliography}
\end{document}